\documentclass[11pt,fleqn]{article}
\usepackage[dvipsnames]{xcolor}
\usepackage{amsmath}
\usepackage{amssymb}
\usepackage{amsfonts}
\usepackage{mathtools}
\usepackage{euscript}
\usepackage{pifont}          
\usepackage{theorem}
\usepackage{charter}
\renewcommand{\leq}{\ensuremath{\leqslant}}
\renewcommand{\geq}{\ensuremath{\geqslant}}
\newcommand{\RR}{\ensuremath{\mathbb R}}
\newcommand{\NN}{\ensuremath{\mathbb N}^*}

\newcommand{\CC}{\ensuremath{\mathbb C}}
\newcommand{\RP}{{\mathbb{R}_+^*}}
\newcommand{\EE}{{\mathsf E}}

\newcommand{\PP}{{\mathsf P}}
\newcommand{\asarrow}{\mbox{$\:\:\stackrel{{\rm a.s.}}%
{\longrightarrow}\:\:$}}

\newcommand{\warrow}{\mbox{$\:\:\Rightarrow\:\:$}}
\newcommand{\conv}{\mbox{$\:\rightarrow\:$}}
\newcommand{\pinf}{\ensuremath{{{+}\infty}}}
\newcommand{\minf}{\ensuremath{{{-}\infty}}}
\definecolor{labelkey}{HTML}{0455BF}
\definecolor{refkey}{rgb}{0,0.6,0.0}
\definecolor{dblue}{HTML}{044EAF}
\definecolor{dgreen}{HTML}{02724A}
\definecolor{myellow}{HTML}{D97904}
\definecolor{dred}{HTML}{D90404}
\usepackage{upref}
\usepackage{hyperref}
\hypersetup{colorlinks=true,%
linktocpage=true,%
linkcolor=dblue,%
citecolor=dgreen,%
urlcolor=dred}

\usepackage{enumitem}

\numberwithin{equation}{section}
\setlist[enumerate]{itemsep=-1pt,topsep=1pt}
\setlist[description]{itemsep=-1pt,topsep=1pt}
\setlist[itemize]{itemsep=-1pt,topsep=1pt}
\usepackage{soul}
\setstcolor{dred}
\usepackage{authblk}
\newcommand{\email}[1]{\href{mailto:#1}{\nolinkurl{#1}}}

\author[1]{Jack W. Silverstein}
\affil[1]{North Carolina State University
\authorcr
Department of Mathematics
\authorcr
Raleigh, NC 27695, USA
\authorcr
~
}
\author[2]{Patrick L. Combettes}
\affil[2]{The City College of The City
University of New York
\authorcr
Department of Electrical Engineering
\authorcr
New York, NY 10031
}
\begin{document}
\title{\sffamily\huge%
Spectral Theory of Large Dimensional Random 
Matrices Applied to Signal Detection\thanks{This technical report
from 1990 is the long version of the paper:\\
J. W. Silverstein and P. L. Combettes, Signal detection via spectral
theory of large dimensional random matrices, 
{\em IEEE Transactions on Signal Processing,}
vol. 40, no. 8, pp. 2100--2105, August 1992,\\
which was published without the proofs. Since several colleagues
have requested those proofs over the years, we make them available
now. Current emails of the authors: 
\email{jack@math.ncsu.edu} and 
\email{plc@math.ncsu.edu}.
}
}

\date{\ttfamily Fall 1990}
\maketitle

\begin{abstract}
Results on the spectral behavior of random matrices as the
dimension increases are applied to the problem of detecting the
number of sources impinging on an array of sensors. A common
strategy to solve this problem is to estimate the multiplicity of
the smallest eigenvalue of the spatial covariance matrix $R$ of the
sensed data from the sample covariance matrix $\widehat{R}$.
Existing approaches, such as that based on information theoretic
criteria, rely on the closeness of the noise eigenvalues of
$\widehat R$ to each other and, therefore, the sample size has to
be quite large when the number of sources is large in order to
obtain a good estimate. The analysis presented in this report
focuses on the splitting of the spectrum of $\widehat{R}$ into
noise and signal eigenvalues. It is shown that, when the number of
sensors is large, the number of signals can be estimated with a
sample size considerably less than that required by previous
approaches. The practical significance of the main result is that
detection can be achieved with a number of samples comparable to
the number of sensors in large dimensional array processing.
\end{abstract}

\newpage

\section{Introduction}
In many signal processing applications, a fundamental problem is
the determination of the number of signals impinging on an array of
sensors. Under the assumption that the vector of sensed data
consists of superimposed random signals corrupted by additive white
noise, the number of signals present in the scene is related to the
multiplicity of the smallest eigenvalue of the spatial covariance
matrix $R$ of the data process, this eigenvalue being equal to the
power of the noise. Since $R$ is unknown, its spectrum must be
approximated by observing that of the sample covariance matrix
$\widehat{R}$ of the data process sampled across time. The
eigenvalues of $\widehat{R}$ being typically distinct\footnote{For
a sufficient condition under which the eigenvalues of $\widehat R$
would be almost surely distinct, see \cite{Okam73}.}, the detection
problem is that of deciding which of the smallest eigenvalues are
associated with the noise. An approach is to use hypothesis tests
on the multiplicity of the smallest eigenvalue of a random matrix,
such as that discussed in \cite{Kris76}. An alternative strategy
based on information theoretic criteria for model selection was
proposed in \cite{Wax85} and was further studied in 
\cite{Kave87,Yin87,Zhan89,Zhao86,Zhao87}. 

All of these detection methods rely on the ergodic theorem and
their performance strongly depends on $R$ being closely
approximated by $\widehat{R}$, requiring the sample size to be
quite large. In applications where the number of signals and,
consequently, the number of sensors, is sizable, the required
number of samples may be prohibitive. The purpose of this report is
to bring into play elements of the spectral theory of random
matrices, more specifically, results on the limiting distribution
of the eigenvalues of random matrices as the dimension increases.
This analysis will show that, when the number of sensors is large,
the number of signals can be estimated with a sample size
considerably less than that required by invoking the ergodic
theorem.

The report is organized as follows. Results from the spectral
theory of random matrices are introduced in Section~\ref{sec:2}.
The application to signal detection is presented in
Section~\ref{sec:3} and numerical results are provided in
Section~\ref{sec:4}. Our concluding remarks appear in
Section~\ref{sec:5}. All of our results are proved in the Appendix
(Section~\ref{sec:6}).

\section{Spectral Theory of Random Matrices}
\label{sec:2}
Throughout this report, $\NN$ will denote the set of strictly
positive integers and $\RP$ the set of strictly positive real
numbers. All the random variables (r.v.'s) are defined on a
probability space $(\Omega,\Sigma,\PP)$. A r.v. $X$ is said to be
in $L^r (\PP)$ ($0\!<\!r\!<\! \pinf$) if $\EE |X|^r < \pinf$. For
r.v.'s, almost sure convergence is denoted by $\asarrow$ and, for
distribution functions\footnote{By a d.f. we mean a
right-continuous nondecreasing function $F$ on $\RR$ with $\lim_{x
\:\rightarrow\:-\infty} F(x) = 0$ and $\lim_{x
\:\rightarrow\:+\infty} F(x) =1$. The support of $F$ is the closed
set ${\mathcal S}=\{x\in\RR \:|\: (\forall \varepsilon \in\RP) 
\:\:\: F(x+\varepsilon)>F(x-\varepsilon) \}$.} (d.f.'s), weak
convergence is denoted by $\warrow$. The transpose of a matrix $A$
is denoted by $A^\top$, its conjugate transpose by $A^*$, and its
trace by ${\rm tr}A$. 

Let $M$ be
an $m\times m$ random matrix with real-valued eigenvalues
$\{\Lambda_{1}, \ldots , \Lambda_{m}\}$. The empirical d.f. of the
r.v.'s $\{\Lambda_{1}, \ldots , \Lambda_{m}\}$ is the stochastic
process defined by\footnote{The characteristic function of a set
$S$ is denoted by $1_S$.} 
\begin{equation}
( \forall \omega \in \Omega)(\forall x \in \RR ) 
\:\:\:\:\: F^{M}(x,\omega) = \frac{1}{m} 
\sum_{i=1}^{m} 1_{\left]- \infty, x \right]} 
(\Lambda_{i} (\omega)).
\end{equation}
We now review the main result, a limit theorem found in
\cite{Yin86}.

{\bfseries Theorem 1 } \cite{Yin86}.
Let $(Y_{ij})_{i,j \geq 1}$ be i.i.d. real-valued r.v.'s with 
$\EE |Y_{11}- \EE Y_{11} |^2=1$. For each $m$ in $\NN$ , 
let $Y_m=[Y_{ij}]_{m \times n}$, where $n=n(m)$ and $m/n\conv y>0$
as $m\conv  \pinf$, and let $T_m$ be an $m\times m$ symmetric
nonnegative definite random matrix independent of the $Y_{ij}$'s
for which there exists a sequence of strictly positive numbers 
$(\mu_k)_{k\in\NN}$ such that for each $k$ in $\NN$ 
\begin{equation}
\label{mom1}
\int_{0}^{+ \infty} x^k dF^{T_m}(x) = \frac{1}{m} {\rm tr} T_m^k
\asarrow \mu_k \:\:\:\: {\rm as} \:\:\:\:m \conv \pinf 
\end{equation} 
and where the $\mu_k$'s satisfy Carleman's
sufficiency condition, $\sum_{k\in\NN} \mu_{2k}^{-1/2k} = \pinf$,
for the existence and the uniqueness of the d.f. $H$ having moments
$(\mu_k)_{k\in\NN}$. Let $M_m=(1/n) Y_m Y_m^\top T_m$. Then, almost
surely, $(F^{M_m})_{m \geq 1}$ converges weakly to a nonrandom d.f.
$F$ having moments
\begin{equation}
(\forall k \in \NN)\:\:\:\nu_k = \sum_{w=1}^{k} y^{k-w} 
\sum \frac{k!}{m_1! \cdots m_w! w!} 
\mu_1^{m_1} \cdots \mu_w^{m_w}
\label{mom}
\end{equation}
where the inner sum extends over all $w$-tuples of positive
integers $(m_1 , \ldots , m_w )$ such that $\sum_{i=1}^{w} m_i =
k-w+1$ and $\sum_{i=1}^{w} im_i = k$. Moreover, these moments
uniquely determine $F$. 

Similar results are given in \cite{Marc67} and \cite{Wach78} with
varying degrees of assumptions, although in both papers the
matrices studied can have complex-valued entries. However, the
proof in \cite{Yin86} can easily be modified to allow
complex-valued entries in $Y_m$ and $T_m$, giving the same result,
provided $T_m$ is Hermitian and we take $M_m=(1/n) Y_mY_m^*T_m$. 

Although it does not appear likely a general explicit expression
for $F$ in terms of $y$ and arbitrary $H$ can be derived, useful
qualitative information can be found from the different methods
used in \cite{Marc67}, \cite{Silv85}, and \cite{Wach78} to express
transforms of $F$ (transforms of Stieltjes type in \cite{Marc67}
and \cite{Wach78}, the characteristic function in \cite{Silv85}).
For example, in \cite{Marc67}, it is shown that the endpoints of
the connected components of the support of $F$ are given by the
extrema of the function 
\begin{equation} 
\label{j7} 
f(\alpha)=-\frac{1}{\alpha}+y\int_{0}^{\pinf}
\frac{dH(x)}{\alpha+1/x}.
\end{equation} 
The analysis in \cite{Silv85} shows how
one can prove that $F$ is absolutely continuous on $\RP$ and
express its derivative, provided the inverse of a certain function
defined by $y$ and $H$ can be analytically extended in the real
part of the complex plane. 

We now provide additional results
apropos of the limiting behavior of $(F^{M_m})_{m \geq 1}$.

{\bfseries Theorem 2.} 
The limiting d.f. $F$ in Theorem 1 is continuous on
$\RP$. Moreover, if $H$ places no mass at $0$ then, almost surely,
$(F^{M_m})_{m \geq 1}$ converges to $F$ uniformly in $\RR$.

{\bfseries Proposition 1.} With the same notation and
hypotheses as in Theorem 1, the following hold: 
\begin{enumerate}
\item  
$F$ and $y$ uniquely determine $H$. 
\item
Almost surely, $(F^{T_m})_{m \geq 1}$ converges to $H$ 
weakly.
\item
$F \warrow H$ as $y\conv 0$.  
\end{enumerate}

Statement (iii) has a direct bearing on the problem of
estimating the spectrum of a covariance matrix from observing that
of a sample covariance matrix. Indeed, the matrix $(1/n)
T_m^{1/2}Y_mY_m^{*}T_m^{1/2}$ (whose eigenvalues are identical to
those of $M_m$\footnote{The reader is reminded that given two
matrices $A_{p \times q}$ and $B_{q \times p}$, where $p\geq q$,
the spectrum of $AB$ is that of $BA$ augmented by $p-q$ zeros.})
encompasses a broad class of sample covariance matrices stemming
from $n$ i.i.d. samples distributed as an $m$-dimensional random
vector $X$ with $\EE X = 0$ and $\EE XX^* = T_m$ (including the
Wishart case when $X$ is multivariate complex Gaussian). In
estimating the spectrum of $T_m$ from the sample covariance matrix,
there seems to be no mention in the literature as to the dependence
of $n$ on $m$, that is, how large the sample size should be
vis-\`{a}-vis the vector dimension in order to estimate the
eigenvalues to within a certain degree of accuracy. Indeed,
asymptotic results are expressed only in terms of the sample size
(see e.g. \cite{Ande84}). The fact that $F$ differs from $H$ for
$y>0$ while $F \warrow H$ as $y\conv 0$, which complements the fact
that, for fixed $m$, $M_m \asarrow T_m$ as $n\conv \pinf$, confirms
the intuitively apparent statement that, for $m$ large, $n$ should
be much larger, in the sense that $m = o(n)$. 

\section{Application to Signal Detection} 
\label{sec:3}

\subsection{Description of the Problem and Assumptions} 

Let $p$ be the number of sensors in the array, $q$ the
unknown number of signals ($q<p$), and $[0,\tau]$ be the
observation interval. At each time $t$ in $[0,\tau]$, the $j$-th
signal present in the scene, the additive noise at the $i$-th
sensor, and the received data at the $i$-th sensor are respectively
represented by the $L^2(\PP)$ complex-valued r.v.'s $S_j(t)$,
$N_i(t)$, and $X_i(t)$. The random vectors $(S(t)=[ S_{1}(t) \ldots
S_{q}(t) ]^{\top})_{t \in [0,\tau]}$ are identically distributed
(i.d.) with nonsingular spatial covariance matrix $R_S = \EE
S(0)S(0)^*$. Moreover, it is assumed that the r.v.'s $(N_i(t) \:|\:
1\! \leq\! i\! \leq\! p,\: t \in [0,\tau] \}$ are independent
and identically distributed (i.i.d.) with $\EE N_1(0)=0$ and $\EE
|N_1(0)|^2=\sigma^2$, where $\sigma^2$ is unknown, and independent
from the r.v.'s $(S_j(t) \:|\:1\!\leq\!j\!\leq\! q,\: t \in
[0,\tau] )$. Let 
\begin{equation}
N(t)= \sigma W(t) = \sigma [ W_{1}(t) \ldots W_{p}(t) ]^{\top} 
\end{equation}
(so that the $W_i(t)$'s are standardized) and
$X(t)=[ X_{1}(t)\ldots X_{p}(t)]^{\top}$. The data collected by the
array of sensors are modeled as observations of the random vector
\begin{equation} 
\label{model} 
X(t) = A S(t) + N(t),\:\:\:\:\:\:\:\:\: t\in [0,\tau], 
\end{equation} 
where $A$ is a $p \times q$
complex matrix depending on the geometry of the array and the
parameters of the signals, and is assumed to have rank $q$. The
detection problem is to estimate $q$ from the observation of $n$
snapshots $\{X(t_1), \ldots , X(t_n)\}$ of the data process. Under
the above assumptions, the random vectors $(X(t))_{t \in [0,\tau]}$
are i.d. with spatial covariance matrix 
\begin{equation} 
R=\EE X(0)X(0)^*=AR_{S}A^{*}+\sigma^{2}I_p,
\end{equation} 
where
$I_p$ denotes the $p \times p$ identity matrix. Moreover, the $p-q$
smallest eigenvalues of $R$ are equal to $\sigma^{2}$. These
eigenvalues will be referred to as the noise eigenvalues and the
remainder of the spectrum will be referred to as the signal
eigenvalues. In practice, $R$ is not known, and its spectrum must
be inferred from observing that of the sample covariance matrix
\begin{equation} 
{\widehat R} = \frac{1}{n} \sum_{i=1}^n
X(t_i)X(t_i)^*. 
\end{equation} 
Loosely speaking, one must then
decide where the observed spectrum splits into noise and signal
eigenvalues. \subsection{General Analysis} For every $t$ in $[0,
\tau ]$, let us assume that the signal vector is given by
\begin{equation} 
\label{signal} 
S(t)=CV(t) \:\:\:\:{\rm with}\:\:\:\:
V(t)=[V_1(t),\ldots,V_q(t)]^\top,
\end{equation} 
where
$C$ is $q\times q$, nonsingular, and the r.v.'s
$\{V_1(t),\ldots,V_q(t)\}$  are i.i.d. with the same d.f. as
$W_1(0)$. It is worth noting that this general formulation
comprises the special case when $S(0)$ is multivariate complex
Gaussian, which is a common assumption in array signal processing.
Let $B=AC$. Then (\ref{model}) yields 
\begin{equation} 
X(t)=\left[ B \:\:\: \sigma I_p \right] 
\left[ \begin{array}{c} V(t)\\ W(t)
\end{array} \right].
\end{equation} 
Notice that $R_S=CC^*$ and
$R=BB^*+\sigma^2I_p$. If we further assume that the $n$ vectors
$\{S(t_1), \ldots , S(t_n)\}$ are independent, then the $n$ data
samples $\{X(t_1), \ldots , X(t_n)\}$ will also be independent and
the corresponding sample covariance matrix $\widehat R$ takes on
the form 
\begin{equation} 
\label{rh} 
\widehat R=\frac{1}{n}[B \:\:\: \sigma
I_p]VV^*[B \:\:\: \sigma I_p]^*,
\end{equation} 
where $V=[V_{ij}]_{(p+q)\times n}$ consists of i.i.d. standardized
entries. 

{\bfseries Theorem 3}. If $W_1(0)$ is standard complex
Gaussian\footnote{A r.v. is said to be standardized complex
Gaussian if its real and imaginary parts are i.i.d. with mean zero
and variance 1/2.}, the joint distribution of the eigenvalues of
$\widehat R$ in (\ref{rh}) is the same as the joint distribution of
the eigenvalues of $\widehat R^{\prime}= (1/n)
Y_pY_p^*(BB^*+\sigma^2I_p)$, where $Y_p$ is any $p\times n$ random
matrix with i.i.d. standardized complex Gaussian entries. In
general, for $p$ and $n$ sufficiently large, with high probability,
the empirical d.f.'s $F^{\widehat R}$ and $F^{\widehat R^{\prime}}$
are close to the d.f. $F$ of Theorem 1 for $m=p$, $y=p/n$, and
$H=F^{BB^*+\sigma^2I_p}$. 

The importance of Theorem~3 becomes immediately apparent. The
observations of the empirical d.f. $F^{\widehat R}$, for suitably
large $p$ and $n$, will not vary very much from one realization to
another, even if $n$ is not large relative to $p$. In fact, by
Theorem 2, with high probability, $F^{\widehat R}$ will be
uniformly close to a d.f. $F$ that depends only on $y$ and the
eigenvalues of $BB^*+\sigma^2I_p$. Hence, a realization of
$F^{\widehat R}$ and the ratio $p/n$ can be used to describe, to
within a certain degree of accuracy, $F^{BB^*+\sigma^2I_p}$, which
will yield $\sigma^2$ and the ratio $y_1=q/p$ which corresponds to
the $q$ strictly positive eigenvalues of $BB^*$. 

Much of the information on the spectrum of
$BB^*+\sigma^2I_p$ can be directly observed from plotting
histograms of the eigenvalues of $\widehat R$, in particular, the
ratio $y_1$ of signal eigenvalues. Let $G$ denote the empirical
d.f. of the eigenvalues of $BB^*+\sigma^2I_p$ which are greater
than $\sigma^2$, and let $b_1$ and $b_2$ denote, respectively, the
smallest and largest of these values. Then, for every $x$ in $\RR$,
we can write 
\begin{equation}
\label{G} 
H(x)=F^{BB^*+\sigma^2I_p}(x)=(1-y_1)1_{[\sigma^2,+
\infty[}(x)+y_1G(x) 
\end{equation} 

{\bfseries Proposition 2.}
When $y<1$, the smallest interval $[x_1,x_4]$ containing the
support of $F$ satisfies $0<x_1<x_4< \pinf$ with $x_1 \uparrow
\sigma^2$ and $x_4 \downarrow b_2$ as $y \downarrow 0$. In
addition, there exists an $\alpha$ in $]-1/\sigma^2,-1/b_1[$ such
that 
\begin{equation} 
\label{split}
g(\alpha) = y\left((1-y_1)\left(\frac
\alpha{\alpha+1/\sigma^2} \right)^2 +
y_1\int_{b_1}^{b_2}\left(\frac
\alpha{\alpha+1/x}\right)^2dG(x)\right)<1 
\end{equation} 
(which can always be found for $y$ sufficiently small) if and only
if the support of $F$ splits into at least two separate components,
with the leftmost interval $[x_1,x_2]$ being a connected component
of the support containing mass $1-y_1$ from $F$. Furthermore, for
$y$ sufficiently small, $x_2 \downarrow \sigma^2$ as 
$y\downarrow 0$ and, if $[x_3,x_4]$ denotes the smallest interval
containing the remaining support of $F$, then $x_3 \uparrow b_1$ as
$y\downarrow 0$. Regardless of the respective location of $x_2$ and
$x_3$ vis-\`{a}-vis $\sigma^2$ and $b_1$, the separation between
the noise and signal portions of the spectrum, i.e. $x_3 - x_2$,
increases as $y$ decreases. When $y>1$, $F$ places mass $1-1/y$ at
the origin, but the remaining support will lie to the right of a
strictly positive value $x_1$. It is still possible for the support
of $F$ to split further provided (\ref{split}) holds. In this case
the leftmost interval $[x_1,x_2]$ will carry mass $(1/y)-y_1$,
leaving mass $y_1$ to the remaining support of $F$ to the right of
$x_2$. When $y=1$ the latter situation applies, except now $x_1=0$,
and there will be no mass at 0.

Thus, if $p$ and $n$ are large enough so that
$F^{\widehat R}$ is close to $F$ with high probability, then for
$y=p/n$ suitably small, an appropriately constructed histogram of
the eigenvalues of $\widehat R$ will display clustering on the left
separated from the rest of the figure. The proportion of the
number of eigenvalues associated with the histogram to the right of
the clustering will then be close to $q/p$, with high
probability.

Although the theory merely guarantees that the
proportion of signal eigenvalues of $\widehat R$ is close to that
of $R$, extensive simulation strongly suggests that the spectrum of
$\widehat R$ splits into two portions containing the exact number
of noise and signal eigenvalues, and that the endpoints of these
portions agree very closely with the ones predicted by the theory.
This point, which will be illustrated in Section 4, leads to the
possibility of the existence of a much stronger underlying spectral
theory deepening the results of Theorem 1. Results along these
lines are known for the extreme eigenvalues when $T_m=\sigma^2I_m$.
Such specific cases will be discussed in Section 3.3. 

Intuitively, the above procedure has advantages over other methods
used to estimate $q$, in particular, those adapted from information
theoretic criteria discussed in \cite{Kave87}, \cite{Wax85},
\cite{Yin87}, \cite{Zhan89}, \cite{Zhao86}, and \cite{Zhao87}. The
latter methods try to exploit the closeness of the noise
eigenvalues of $\widehat R$ to each other as well as their
separation from the remaining signal eigenvalues. Usually the
sample size has to be quite large for the smaller eigenvalues to
cluster. On the other hand, only the separation of the two classes
of eigenvalues is needed when viewing the spectrum, so a suitable
$n$ can conceivably be much smaller, sometimes even smaller than
$p$. In other words, previous methods require $\widehat R$ to be
near $BB^*+\sigma^2I_p$, while, for situations where $p$ is
sizable, the present analysis requires $n$ to be large enough so
that the support of $F$ separates. 

\subsection{Specific Cases} 
An
important case to consider is the one for which no signal is
present, that is, when $B=0$, or equivalently, when
$T_m=\sigma^2I_m$. Then it is known
\cite{Silv77,Jons82,Marc67} that, for $y\leq 1$, $F$
is continuously differentiable, where 
\begin{equation} 
\label{dens} 
F'(x)=
\left\{ \begin{array}{ll} {{ \text{$ \dfrac{\left((x -
\sigma^2(1\!-\!\sqrt y )^2)(\sigma^2(1\!+\! \sqrt y )^2-x)
\right)^{1/2} } {2\pi\sigma^2yx}$}}} &\text{if
$\sigma^2(1\!-\!\sqrt y)^2<x<\sigma^2(1\!+\!\sqrt y)^2$};\\ 
0 &\text{otherwise} 
\end{array} \right. 
\end{equation}
and for $y>1$, $F$ has derivative (\ref{dens}) on $\RP$ and mass
$1-1/y$ at 0. Furthermore, the largest eigenvalue of $M_m$
converges almost surely [respect. in probability] to
$\sigma^2(1+\sqrt y)^2$ as $m\conv \pinf$ if and only if $\EE
Y_{11}=0$ and $Y_{11} \in L^4(\PP)$ [respect. $x^4 \PP\{
\omega\in\Omega \:|\:\:|Y_{11}(\omega)|\ge x \}\conv 0$ as $x\conv
\pinf]$ \cite{Bai88,Gema80,Silv90,Yin88}. The
almost sure convergence of the smallest eigenvalue of $M_m$ to
$\sigma^2(1-\sqrt y)^2$ when $y<1$ has thus far been shown only for
$Y_{11}$ standardized Gaussian \cite{Silv88} (it is remarked here
that the results on the extreme eigenvalues have been verified for
$Y_{11}$ real-valued, but, again, the proofs can be extended to the
complex case). 

The above results can be used to investigate the
possibility of no signals arriving at the sensors. Certainly, the
existence of at least one signal would be in doubt if the number of
samples were quite large but histograms indicate only one connected
component away from 0. But, for any $y$, comparisons can be made
between histograms of the eigenvalues, (\ref{dens}), and $F'$ when
$B\neq0$, to infer whether or not signals are present, provided the
latter densities exist and appear different enough from (\ref{dens})
to make a distinction. 

For this reason it is mentioned briefly
here the case when $G$ places mass at one value, $b>\sigma^2$.
Except for the situation of only one signal, this case is not
typically found in practice. However, the d.f. $F$ can be
completely determined and its properties strongly suggest the
smoothness and appearance of $F$ for general $G$. Only the case
$y<1$ will be outlined, the remaining cases for $y$ following as
above. From the analysis in \cite{Silv85} it can be shown that $F$
is continuously differentiable with derivative of the form
\begin{equation} 
\label{otherw} 
F^{\prime}(x)= \left\{ \begin{array}{ll} k {\Large
{\dfrac{ \left( p_3(x)+x \sqrt {p_4(x)} \right)^{1/3}-
\left( p_3(x)-x \sqrt{p_4(x)} \right)^{1/3}  }x}} &\text{if}\;\;
p_4(x)\geq 0;\\[5mm] 
0 &\text{otherwise.} 
\end{array}
\right. 
\end{equation} 
Here, $p_3$ and $p_4$ are, respectively,
third and fourth degree polynomials depending continuously on $y$,
$y_1$, $\sigma^2$, $b$, and the leading coefficient of $p_4$ is
negative.  The latter polynomial has either two real roots,
$0<x_1<x_4$, so that $F$ has support on $[x_1,x_4]$, or four real
roots, $0<x_1<x_2<x_3<x_4$, which is the above mentioned case where
the support of $F$ splits into two intervals, $[x_1,x_2]$ and
$[x_3,x_4]$, with $F(x_2)-F(x_1)=1-y_1$. Using (\ref{split}), it is
straightforward to show $F$ splits if and only if 
\begin{equation}
\label{third}
y\frac{\left((b^2y_1)^{1/3}+(\sigma^4
(1-y_1))^{1/3}\right)^3}{(b-\sigma^2)^2} <1.
\end{equation} 
When the left-hand side of (\ref{third}) is equal to 1,
then $p_4$ still has four real roots, but $x_2=x_3$.  
When
(\ref{third}) holds, $F'$ is unimodal on each of the intervals,
with infinite slopes at each endpoint. If there is a $y<1$, say
$y_o$, for which the left-hand side of (\ref{third}) is equal to 1,
then, since the graph of $F'$ varies continuously with $y$, as $y$
increases from 0, the separate curves eventually join (at $y=y_o$)
and the single curve will display two relative maxima, at least for
$y$ near $y_o$. Thus, although $y$ may not be small enough to
split $F'$, it may still be possible to infer the number of signal
eigenvalues from the shape of a histogram. 

\section{Simulation Results} 
\label{sec:4}
The objective of this section is to illustrate some
aspects of our analysis through their application to the case of a
linear array with $p$ sensors receiving noisy signals from $q$
narrow-band far-field sources. The sensors are assumed to be
omnidirectional with unity gain and uniform spacing $\lambda/2$,
where $\lambda$ is the signal wavelength\footnote{In this context,
the matrix $A$ in (\ref{model}) has a Vandermonde structure with
$A_{ki}=\exp(-\imath\pi(k-1)\sin\theta_i)$ ($1\!\leq\!k\!\leq\!p$,
$1\!\leq\!i\!\leq\!q$), where $\theta_i$ is the angle of arrival of
the $i$-th signal with respect to the normal to the array.}. 

Our analysis applies to cases where $p$ is large. Simulations have
supported its applicability for values of $p$ as low as 30. In the
simulation presented here, the number of sensors is set to $p=50$
and the noise is zero mean, white, complex Gaussian, with power
$\sigma^2 = 1$. The signal scenario consists of $q=35$ partially
correlated sources with angles of arrivals uniformly spaced between
$-70^o$ and $70^o$ and power selected at random from a uniform
distribution so as to yield signal-to-noise ratios ranging from 0dB
to 10dB. The signal vector is multivariate complex Gaussian and
obtained according to (\ref{signal}), where $C$ is a randomly
generated banded matrix.

The spectrum ${\mathcal L}$ of $R = BB^*+I_{50}$, where $B=AC$, was
computed in order to obtain an explicit
expression for the functions $f(\cdot)$ of (\ref{j7}) and
$g(\cdot)$ of (\ref{split}). Newton's method was used to find the
minimum of $g(\cdot)$ over $]-1/\sigma^2,-1/b_1[$ and, whence, it
was found that the largest value of $y$ for which the splitting of
the spectrum occurs (i.e. (\ref{split}) holds) is $\tilde{y}=
1.058$. Then, with the above configuration, four experiments were
performed with the following number of samples $n$: 50, 100, 250,
and 1500 (which corresponds to values of $y$ of 1, 1/2, 1/5, and
1/30, respectively). In each experiment, 10 realizations ${\mathcal
L}_1, \ldots, {\mathcal L}_{10}$ of the spectrum of the sample
covariance matrix $\widehat R$ were observed, the eigenvalues being
arranged in nondecreasing order. The results of these experiments
are shown in Tables 1 through 4. Even for $y=1$, the $p-q=15$
smallest eigenvalues are seen to cluster to the left of most of the
observed spectra. This confinement delimitates exactly the noise
portion of the spectrum and, thereby, detects the exact number of
signals. As discussed earlier, for a given value of $y$, the
theoretical endpoints of the supports of the noise and signal
portions of the spectrum can be determined from the location of the
relative extrema of $f(\cdot)$. Newton's method was used to this
end and gave the results shown in Table 5. In agreement with
Proposition 2, it is seen that, as $y$ decreases, the separation
$x_3-x_2$ increases while the endpoints converge towards the
theoretical values. 

\begin{center}
{\bfseries Table 1. Observed Spectra - $y=1$. \\}
\vspace{0.4cm}
\begin{tabular}{|c||c|c|c|c|c|c|c|c|c|c||c|} \hline
 &${\mathcal L}_1$&${\mathcal L}_2$&${\mathcal L}_3$&${\mathcal
L}_4$&${\mathcal L}_5$&${\mathcal
 L}_6$&${\mathcal L}_7$&${\mathcal L}_8$&${\mathcal
L}_9$&${\mathcal L}_{10}$&${\mathcal L}$\\
 \hline\hline$\lambda_{1}$
 &0.00&0.00&0.00&0.00&0.00&0.00&0.00&0.00&0.00&0.00&1\\
 $\lambda_{2}$  &0.00&0.01&0.01&0.01&0.01&0.01&0.01&0.00&0.00&0.00&1\\
 $\vdots$       &&&\vdots&&&&&\vdots&&&\vdots\\
 $\lambda_{10}$ &0.31&0.33&0.32&0.32&0.28&0.28&0.33&0.35&0.34&0.34&1\\
 $\lambda_{11}$ &0.43&0.40&0.37&0.41&0.38&0.36&0.38&0.40&0.45&0.50&1\\
 $\lambda_{12}$ &0.45&0.47&0.48&0.49&0.44&0.49&0.41&0.52&0.48&0.55&1\\
 $\lambda_{13}$ &0.57&0.57&0.50&0.60&0.58&0.61&0.58&0.63&0.64&0.64&1\\
 $\lambda_{14}$ &0.67&0.64&0.64&0.74&0.80&0.73&0.78&0.75&0.74&0.73&1\\
 $\lambda_{15}$ &0.86&0.87&0.83&1.05&0.96&0.95&0.95&0.86&0.87&0.90&1\\
 $\lambda_{16}$ &1.38&1.64&1.40&1.90&1.18&1.45&2.35&1.71&1.61&1.70&5.34\\
 $\lambda_{17}$ &2.59&2.72&2.41&2.81&1.82&3.41&3.20&2.50&2.85&2.40&6.20\\
 $\lambda_{18}$ &5.61&5.21&4.74&4.97&3.52&4.11&5.85&4.88&4.67&5.32&21.4\\
 $\lambda_{19}$ &7.98&7.64&8.22&7.37&9.66&6.06&7.03&8.47&6.70&6.07&23.1\\
 $\lambda_{20}$ &11.4&9.87&10.8&9.67&11.6&8.16&11.1&11.1&12.2&11.8&25.7\\
 $\lambda_{21}$ &14.8&13.3&11.9&11.7&16.6&13.7&14.3&12.9&14.5&13.3&49.2\\
 $\vdots$       &&&\vdots&&&&&\vdots&&&\vdots\\
 $\lambda_{49}$ &1159&1074&1137&1065&1067&1154&1128&1123&1135&1229&756\\
 $\lambda_{50}$ &1470&1309&1233&1390&1458&1556&1547&1306&1347&1522&932\\
 \hline
 \end{tabular}
 \end{center}
 \newpage
 \begin{center}
 {\bfseries Table 2. Observed Spectra - $y=1/2$. \\}
 \vspace{0.4cm}
 \begin{tabular}{|c||c|c|c|c|c|c|c|c|c|c||c|} \hline
  &${\mathcal L}_1$&${\mathcal L}_2$&${\mathcal L}_3$&${\mathcal
L}_4$&${\mathcal L}_5$&${\mathcal
  L}_6$&${\mathcal L}_7$&${\mathcal L}_8$&${\mathcal
L}_9$&${\mathcal L}_{10}$&${\mathcal L}$\\
  \hline\hline$\lambda_{1}$
  &0.22&0.23&0.23&0.23&0.22&0.21&0.21&0.22&0.21&0.22&1\\
  $\lambda_{2}$  &0.28&0.29&0.25&0.25&0.28&0.25&0.29&0.26&0.26&0.26&1\\
  $\vdots$       &&&\vdots&&&&&\vdots&&&\vdots\\
  $\lambda_{10}$ &0.77&0.79&0.75&0.71&0.79&0.71&0.75&0.73&0.76&0.70&1\\
  $\lambda_{11}$ &0.97&0.91&0.83&0.77&0.84&0.79&0.78&0.84&0.82&0.79&1\\
  $\lambda_{12}$ &1.05&0.96&0.90&0.91&0.89&0.88&0.89&0.93&0.90&0.89&1\\
  $\lambda_{13}$ &1.11&1.06&0.96&0.99&0.99&0.95&0.96&0.94&0.93&0.97&1\\
  $\lambda_{14}$ &1.20&1.16&1.07&1.10&1.09&1.11&1.14&1.13&1.01&1.03&1\\
  $\lambda_{15}$ &1.41&1.32&1.20&1.30&1.16&1.24&1.31&1.31&1.21&1.33&1\\
  $\lambda_{16}$ &3.21&3.21&3.35&3.81&3.01&3.14&3.35&3.29&2.97&3.54&5.34\\
  $\lambda_{17}$ &4.04&4.27&4.49&4.78&4.07&4.25&4.82&5.36&3.65&4.61&6.20\\
  $\lambda_{18}$ &13.2&10.2&11.8&11.7&13.4&11.8&12.1&11.5&11.2&11.6&21.4\\
  $\lambda_{19}$ &13.7&14.7&15.9&16.5&14.4&13.7&15.3&16.8&13.4&17.7&23.1\\
  $\lambda_{20}$ &16.4&18.1&17.5&21.1&18.4&19.2&19.1&19.5&15.7&19.8&25.7\\
  $\lambda_{21}$ &29.4&30.1&28.1&32.7&29.8&25.9&32.0&31.8&28.1&28.7&49.2\\
  $\vdots$       &&&\vdots&&&&&\vdots&&&\vdots\\
  $\lambda_{49}$ &1019& 959& 879&1137&1012& 941&1083& 894& 910& 946&756\\
  $\lambda_{50}$ &1457&1174&1030&1360&1117&1053&1137&1186&1149&1084&932\\
  \hline
  \end{tabular}
  \end{center}

  \begin{center}
  {\bfseries Table 3. Observed Spectra - $y=1/5$. \\}
\vspace{0.4cm}
  \begin{tabular}{|c||c|c|c|c|c|c|c|c|c|c||c|} \hline
   &${\mathcal L}_1$&${\mathcal L}_2$&${\mathcal L}_3$&${\mathcal
L}_4$&${\mathcal L}_5$&${\mathcal
   L}_6$&${\mathcal L}_7$&${\mathcal L}_8$&${\mathcal
L}_9$&${\mathcal L}_{10}$&${\mathcal L}$\\
   \hline\hline$\lambda_{1}$
   &0.54&0.52&0.48&0.53&0.51&0.52&0.52&0.47&0.50&0.49&1\\
   $\lambda_{2}$  &0.60&0.56&0.55&0.57&0.57&0.58&0.54&0.59&0.57&0.54&1\\
   $\vdots$       &&&\vdots&&&&&\vdots&&&\vdots\\
   $\lambda_{10}$ &0.94&0.93&0.94&0.90&0.92&0.95&0.94&0.96&0.91&0.92&1\\
   $\lambda_{11}$ &1.01&0.98&1.00&0.97&0.99&0.97&0.95&1.06&1.00&0.97&1\\
   $\lambda_{12}$ &1.09&1.03&1.06&1.02&1.02&1.04&1.01&1.10&1.09&1.01&1\\
   $\lambda_{13}$ &1.12&1.06&1.13&1.11&1.07&1.10&1.11&1.14&1.12&1.10&1\\
   $\lambda_{14}$ &1.18&1.12&1.24&1.16&1.17&1.17&1.19&1.18&1.25&1.17&1\\
   $\lambda_{15}$ &1.36&1.27&1.31&1.32&1.21&1.38&1.26&1.27&1.34&1.19&1\\
   $\lambda_{16}$ &4.35&4.78&4.37&4.84&4.96&4.03&4.26&4.19&4.71&4.36&5.34\\
   $\lambda_{17}$ &5.26&6.05&6.31&6.13&5.64&5.54&5.61&5.71&5.29&4.94&6.20\\
   $\lambda_{18}$ &17.4&17.5&18.1&18.5&17.0&17.4&16.8&17.2&17.7&16.5&21.4\\
   $\lambda_{19}$ &19.4&20.1&18.7&20.4&18.9&20.0&19.0&20.3&19.2&18.8&23.1\\
   $\lambda_{20}$ &22.9&25.1&21.3&22.3&22.7&21.9&23.4&22.4&22.9&23.3&25.7\\
   $\lambda_{21}$ &36.7&40.2&39.6&39.8&40.8&39.8&42.4&40.9&41.9&37.0&49.2\\
   $\vdots$       &&&\vdots&&&&&\vdots&&&\vdots\\
   $\lambda_{49}$ & 902& 856& 887& 875& 800& 818& 893& 889& 878& 831&756\\
   $\lambda_{50}$ &1063& 985&1004&1064&1142& 980&1043&1068&1008& 993&932\\
   \hline
   \end{tabular}
   \end{center}
   \newpage
   \begin{center}
   {\bfseries Table 4. Observed Spectra - $y=1/30$. \\}
\vspace{0.4cm}
\begin{tabular}{|c||c|c|c|c|c|c|c|c|c|c||c|} \hline
 &${\mathcal L}_1$&${\mathcal L}_2$&${\mathcal L}_3$&${\mathcal
L}_4$&${\mathcal L}_5$&${\mathcal
 L}_6$&${\mathcal L}_7$&${\mathcal L}_8$&${\mathcal
L}_9$&${\mathcal L}_{10}$&${\mathcal L}$\\
\hline\hline$\lambda_{1}$
&0.82&0.82&0.81&0.82&0.83&0.81&0.83&0.82&0.81&0.83&1\\
$\lambda_{2}$  &0.84&0.84&0.84&0.86&0.85&0.83&0.85&0.85&0.85&0.85&1\\
$\vdots$       &&&\vdots&&&&&\vdots&&&\vdots\\
$\lambda_{10}$ &1.02&1.02&1.01&1.01&1.01&1.00&1.02&1.02&1.02&1.02&1\\
$\lambda_{11}$ &1.04&1.04&1.03&1.03&1.04&1.03&1.04&1.04&1.04&1.06&1\\
$\lambda_{12}$ &1.06&1.04&1.06&1.06&1.07&1.05&1.06&1.06&1.05&1.07&1\\
$\lambda_{13}$ &1.08&1.07&1.09&1.08&1.09&1.08&1.09&1.07&1.08&1.09&1\\
$\lambda_{14}$ &1.10&1.11&1.11&1.10&1.12&1.11&1.10&1.08&1.10&1.11&1\\
$\lambda_{15}$ &1.16&1.13&1.13&1.14&1.17&1.14&1.14&1.11&1.15&1.13&1\\
$\lambda_{16}$ &5.32&5.09&5.23&5.27&5.18&5.12&5.03&5.27&5.12&5.52&5.34\\
$\lambda_{17}$ &5.97&6.34&6.08&6.10&5.94&6.16&6.26&5.96&6.05&5.92&6.20\\
$\lambda_{18}$ &20.5&20.6&20.1&19.4&20.2&19.8&21.7&20.5&20.1&20.9&21.4\\
$\lambda_{19}$ &22.0&22.1&21.7&21.3&22.4&22.4&23.1&23.0&22.6&22.2&23.1\\
$\lambda_{20}$ &24.8&25.3&26.0&25.5&25.2&24.2&25.5&25.3&25.1&24.5&25.7\\
$\lambda_{21}$ &48.6&48.8&49.6&47.7&48.3&48.9&46.8&47.9&47.7&48.4&49.2\\
$\vdots$       &&&\vdots&&&&&\vdots&&&\vdots\\
$\lambda_{49}$ & 764& 807& 818& 766& 772& 765& 779& 797& 766& 788&756\\
$\lambda_{50}$ & 944& 978& 929& 947& 896& 948& 966& 956& 991& 962&932\\
\hline
\end{tabular}

\vskip 5mm

{\bfseries Table 5. Theoretical Bounds for Noise and Signal
Spectrum Supports. \\}
\vspace{0.2cm}
\begin{tabular}{|c||c|c|c||c|} \hline
         & $y=1$ & $y=1/5$ & $y=1/30$ & $y=0$ \\ \hline\hline
$x_{1}$  & 0.000 & 0.4642  & 0.789    & 1.000 \\
$x_{2}$  & 1.124 & 1.369   & 1.184    & 1.000 \\
$x_{3}$  & 1.167 & 3.970   & 5.785    & 5.342 \\
$x_{4}$  & 1586  & 1137    & 995.4    & 931.6 \\ \hline
\end{tabular}
\end{center}

\section{Conclusion} 
\label{sec:5}
We have applied results
from the spectral theory of large dimensional random matrices to
the signal detection problem in situations where the number of
sources is sizable. A theoretical foundation was established for
the analysis of the splitting of the spectrum of the sample
covariance matrix between a connected noise component and a
remaining signal component. While conventional methods require that
the sample size be impracticably large in order to closely
approximate the spatial covariance matrix, the present analysis
shows that the observed spectrum will split with high probability
with a number of samples comparable to the number of sensors. As
far as the detection problem is concerned, the eigenvalues of the
spatial covariance matrix $R$ need not be estimated with a high
degree of precision; only the accurate splitting of the spectrum is
required.

This work should suggest to the engineering community that by
simply observing the spectrum of a large dimensional sample
covariance matrix, highly relevant information can be extracted
when the sample size is not exceedingly large. In the context of
large dimensional array processing, the practical significance of
our main result is that detection can be achieved when the sample
size is only on the same order of magnitude as the number of
sensors.

\section{Appendix: Proofs}
\label{sec:6}
The imaginary part of a complex number $z$ is denoted by $\Im z$.

{\bfseries Proof of Proposition 1.} 
(i) is established by noting that the sequence $(\nu_k)_{k\in\NN}$
can be derived from $F$ and, from (\ref{mom}), the sequence
$(\mu_k)_{k\in\NN}$, which uniquely determines $H$, can be
computed unambiguously; (ii) follows from (\ref{mom1}) and the
Fr\'{e}chet-Shohat theorem\footnote{The Fr\'{e}chet-Shohat theorem
\cite{Loev77} states that if $( F_n)_{n \geq 1}$ is a sequence of
d.f.'s having moments of all orders with $(\forall k\in\NN) \:\:$
$\lim_{n \:\rightarrow\:+\infty} \int_{-\infty}^{+\infty} x^k
dF_n(x) = \mu_k$ finite, and if $F$ is the only d.f. with moments
$( \mu_k)_{k\in\NN}$, then $F_n \warrow F$ as $n \conv \pinf$.};
(iii) is a direct consequence of the Fr\'{e}chet-Shohat theorem
since (\ref{mom}) implies that 
$(\forall k\in\NN)$ $\nu_k \conv \mu_k$ as $y \conv 0$. 

\begin{figure}
\begin{center}
\includegraphics[width=14.0cm]{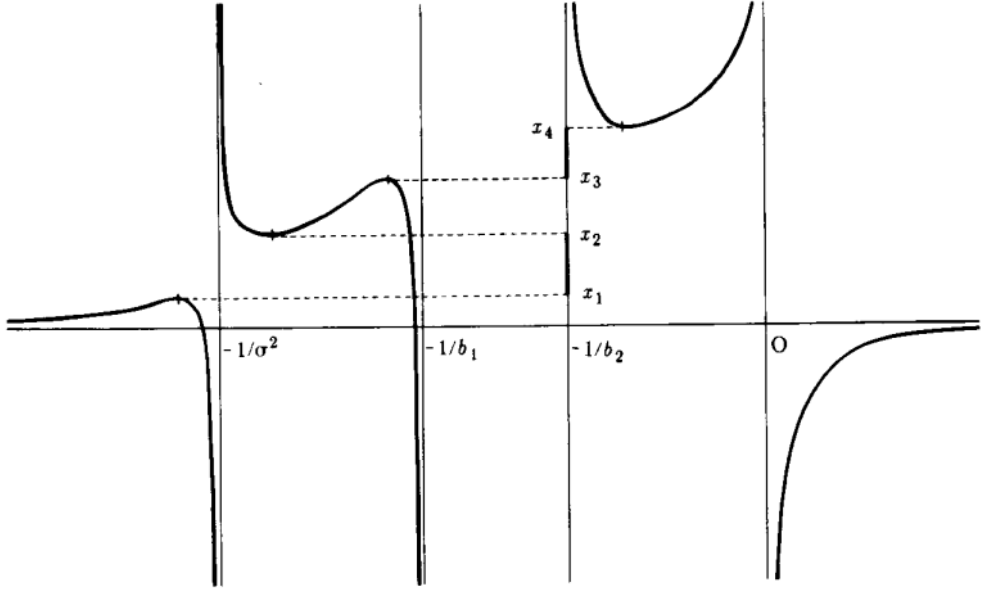}
\end{center}
\caption{ {Graph of $f(\cdot)$.}}
\label{fig:1}
\end{figure}

{\bfseries Proof of Proposition 2.} 
In \cite{Marc67} the matrix corresponding to $M_m=(1/n)X_mX_m^*T_m$ is
$B_n=(1/n)X_m^*T_mX_m$. The spectra of these two matrices differ
only by an additional number of zeros, $n-m$ extra zeros to the
spectrum of $B_n$ when $m<n$, $m-n$ to the spectrum of $M_m$ when
$m>n$. Let $K$ denote the limiting empirical d.f. of the eigenvalues of
$B_n$. It follows that
\begin{equation}
\label{j1}
(\forall x\in\RR)\:\:\:\:K(x)=(1-y)1_{[0,+\infty[}(x)+yF(x).
\end{equation}
Let ${\mathcal A} (\cdot)$ be the Stieltjes transform of $K$, i.e.
\begin{equation}
\label{j2}
(\forall z \in \CC\smallsetminus\RR ) \:\:\:\: {\mathcal A} (z)=
\int_{-\infty}^{+ \infty} \frac{dK(x)}{x-z},
\end{equation}
and let $ {\mathcal B}(\cdot)$ denote the Stieltjes transform of $F$.
From (\ref{j1}) and (\ref{j2}), it follows that
\begin{equation}
\label{j3}
{\mathcal A} (z)=\frac{1-y}{-z}+y {\mathcal B}(z). 
\end{equation}
In \cite{Marc67} it is shown that ${\mathcal A} (z)$, for $\Im\,z>0$,
is the unique solution to the equation
\begin{equation}
\label{j4}
z=-\frac1{{\mathcal A} (z)}+y
\int_0^{+ \infty}\frac{x dH(x)}{1+x {\mathcal A} (z)}, 
\end{equation}
from which $K$ can be calculated from
\begin{equation}
\label{j5}
K(x_2)-K(x_1)=\lim_{\eta \downarrow 0}\frac1{\pi}\int_{x_1}^{x_2}
\Im {\mathcal A} (\xi+\imath\eta)d\xi,
\end{equation}
where $x_1$ and $x_2$ are continuity points of $K$.   Associated
with the above inversion formula is the following.  If $w_1<w_2$
are values lying outside the support of $K$, then
\begin{equation}
\label{j6}
K(w_2)-K(w_1)=-\frac1{2\pi \imath}\oint_C {\mathcal A} (z)dz,
\end{equation}
where $C$ can be taken as the circle in the complex plane having a
diameter with endpoints $w_1$ and $w_2$ on the real axis. 
It is remarked in \cite{Marc67} that, on the union of intervals on
the real axis outside the support of $K$, ${\mathcal A} $ is real and
strictly increasing, and is continuous on each interval. Therefore,
its inverse exists on the range of these intervals and is given by
(\ref{j4}) for ${\mathcal A}$ and $z$ real. This inverse, denoted by
$f(\cdot)$, is given by (\ref{j7}). Writing $H$ as in (\ref{G}),
$f(\cdot)$ takes on the form
\begin{equation}
\label{j17}
f(\alpha)=-\frac{1}{\alpha} + \frac{y(1-y_1)}{\alpha+ 1/\sigma^2}+yy_1
\int_{b_1}^{b_2}\frac{dG(x)}{\alpha + 1/x}.
\end{equation} 
The asymptotic properties of $x_1$ through $x_4$, the existence of
$x_2$ and $x_3$, when $g(\alpha)$ in (\ref{split}) is less than 1,
and the existence of mass at $0$ when $y>1$, all follow from
elementary calculus\footnote{In the case when the integral
$I(\alpha)$ in (\ref{j17}) satisfies $I(b_1-)= \minf$ and $I(b_2+)=
\pinf$, the results are straightforward.  If $I(b_1-)$ [respect.
$I(b_2+)$] is finite, the additional fact that $f(\cdot)$ does not
exist in any interval $]b_1,b_1+\varepsilon[$ [respect.
$\:]b_2-\varepsilon,b_2[ \:\:$] is needed. The latter is proven by
using the fact that whenever $I(\alpha)$ exists,
$G^{\prime}(\alpha) = 0$, which can be verified in a
straightforward manner. Thus, if $f(\cdot)$ were to exist on, say
$]b_1,b_1+\varepsilon[$, then necessarily $G$ would place mass at
$b_1$, resulting in $I(b_1-)=\minf$, a contradiction.}, together
with the fact that $f^{\prime}(\alpha)=(1-g(\alpha))/\alpha^2$.
Figure~\ref{fig:1} shows a typical graph of $f$ when separation
occurs. 

As for the mass $F$
assigns to this interval and to the remaining portions of the support,
(\ref{j6}) can be used together with a change of variables from
(\ref{j4}).  
We will only derive the mass for $[x_1,x_2]$ when $y<1$, the other
portions of the support and cases ($y>1$, $y=1$) being similar.  
With $0<w_1<x_1$, and $w_2$ lying
slightly to the right of $x_2$ we have
\begin{equation}
K(w_2)-K(w_1)=-\frac1{2\pi \imath}\oint_C \left(\frac1{\alpha}-
\frac{\alpha y(1-y_1)} {(\alpha+1/\sigma^2)^2}-yy_1\int_{b_1}^{b_2}
\frac{\alpha dG(x)} {(\alpha+ 1/x)^2} \right) d \alpha,
\end{equation}
where $C$ is a simple closed positively oriented contour enclosing
$-1/\sigma^2$ but not the origin nor any part of $[-1/b_1,-1/b_2]$.
From Cauchy's integral formula, we have 
\begin{equation}
K(w_2)-K(w_1)=y(1-y_1), 
\end{equation}
so that from (\ref{j1}) we conclude $[x_1,x_2]$ contains mass
$1-y_1$ from $F$. Let us mention here that the values of $F(0)$ are
obtained in the proof of Theorem 2.  

The following lemma will be needed in the proof of Theorem 2.

{\bfseries Lemma 1.}
Let $A$ and $B$ be Hermitian, nonnegative definite matrices. Then
for every $\alpha$ and $\beta$ in $\RP$, $F^{AB}(\alpha\beta) 
\leq F^A(\alpha)+F^B(\beta)$. 

{\bfseries Proof.} 
For any $m \times m$ matrix $C$ with real eigenvalues,
let $\ell^C(i)$ be the $i$-th largest eigenvalue of $C$ if 
$1\leq i\leq m$, 
and 0 otherwise. From a routine extension of a result in
\cite{Fan51}, for every positive integers $r_A$ and $r_B$ 
\begin{equation}
\ell^{AB}(m-(r_A+r_B)) \geq \ell^A(m-r_A) \ell^B(m-r_B).
\end{equation}
In particular, if we let $r_A$ [respect. $r_B$] be the number of
eigenvalues of $A$ [respect. $B$] less than or equal to $\alpha$
[respect. $\beta$], the result follows. 

{\bfseries Proof of Theorem 2.} 
We shall use here an argument similar to that used in
\cite{Wach78}. The fact that $K$, and therefore $F$, is continuous
on $\RP$ can be proven by contradiction.  First we notice from
(\ref{j4}) that ${\mathcal A}(z)$ satisfies
\begin{equation}
\label{j9}
{\mathcal A}(z)= \left( -z+y\int_0^{+ \infty}\frac{x dH(x)}{1+x 
{\mathcal A}(z)} \right)^{-1}. 
\end{equation}
Suppose $x_0>0$ is a discontinuity point of $K$ with jump $\mu$. Then
(\ref{j2}) gives 
\begin{equation}
\Im {\mathcal A}(x_0 \!+ \!\imath\eta)=
\int_{- \infty}^{+\infty} \frac{\eta dK(x)} {(x \!-
\!x_0)^2+\eta^2} \geq
\int_{\{x_0\}} \frac{ \eta dK(x)} {(x \!- \!x_0)^2+\eta^2} =
\frac{\mu}{\eta}\conv + \infty \:\:\: {\rm as} \:\:\:
\eta\downarrow 0.
\end{equation}
On the other hand,
\begin{equation}
\left|\int_0^{\pinf}
\frac{x dH(x)}{1+x {\mathcal A}(x_0+\imath\eta)}
\right|\leq\frac1{\Im {\mathcal A}(x_0+\imath\eta)}\conv 0 \:\:\:
{\rm as} \:\:\: \eta \downarrow 0.
\end{equation}
The contradiction then arises from the fact that as $\eta
\downarrow 0$ the left-hand side of (\ref{j9}) becomes unbounded,
while the right-hand side approaches $-1/x_0$. 
The last assertion
is proved by noting that if $H$ places no mass at $0$, $(F^{M_m}(0
\pm))_{m \geq 1}$ converges almost surely to $F(0 \pm)$. Indeed,
trivially, $(\forall m \in \NN)\:\:F^{M_m}(0-)=F(0-)=0$. Moreover,
since the eigenvalues of $M_n$, $(1/n)Y_mY_m^*$, and $T_m$ are all
positive, it follows from Sylvester's inequality \cite{Gant77}
that 
\begin{equation}
\label{j11}
\max \{ F^{(1/n)Y_mY_m^*}(0),F^{T_m}(0)\} \leq F^{M_n}(0) \leq 
F^{(1/n)Y_mY_m^*}(0)+F^{T_m}(0).
\end{equation}
Since 0 is not a mass point of $H$, it is a point a continuity of
$H$ and, therefore, from (ii) in Proposition 1, almost surely,
$(F^{T_m}(0))_{m \geq 1}$ converges to $0$. Thus, from (\ref{j11}), 
\begin{equation}
\label{j15}
|F^{(1/n)Y_mY_m^*}(0)-F^{M_m}(0)| \asarrow 0 
\end{equation}
When $y \leq 1$, from (\ref{dens}), $F^{(1/n)Y_mY_m^*}(0) \asarrow
0$. Fix an arbitrary $\delta$ in $\RP$. Let $\varepsilon$ be in
$]0,(1-\sqrt{y})^2[$ such that $F^H(\varepsilon) < \delta$. By
(\ref{dens}) again, $F^{(1/n)Y_mY_m^*}(\varepsilon) \asarrow 0$.
Moreover, from (ii) in Proposition 1, almost surely,
$\limsup_{m \rightarrow\pinf} F^{T_m}(\varepsilon) < \delta$. Then, by 
Lemma 1, the almost sure limit of $F^{M_n}(\varepsilon^2)$ is less
than $\delta$ and it follows that $F(0)<\delta$. Since $\delta$ can
be made arbitrarily small, we conclude that
\begin{equation}
F^{M_m}(0) \asarrow 0=F(0). 
\end{equation}
When $y>1$, from Section 4.2, $F^{(1/n)Y_mY_m^*}(0) \asarrow
1-(1/y)$. Therefore, it will follow from (\ref{j15}) and the same
argument that 
\begin{equation}
F^{M_m}(0) \asarrow 1-\dfrac{1}{y}=F(0). 
\end{equation}
The proof is
complete since, if a sequence of d.f.'s $( F_n)_{n \geq 1}$
converges weakly to a d.f. $F$ and if $(F_n(x\pm ))_{n \geq 1}$
converges to $F(x \pm)$ at every point $x$ of discontinuity of $F$,
then $(F_n)_{n \geq 1}$ converges to $F$ uniformly in $\RR$
\cite{Chow88}. 

{\bfseries Proof of Theorem 3}. First, notice
that the eigenvalues of $\widehat R$ are the same as the $p$
largest eigenvalues of 
\begin{equation} 
\label{eig5} 
\frac{1}{n}VV^*[B \:\:\:
\sigma I_p]^*[B \:\:\: \sigma I_p].
\end{equation} 
Thus,
for $p$ and $n$ sufficiently large, we see that the empirical d.f.
$F^{\widehat R}$ of the eigenvalues of $\widehat R$, together with
$q$ zeros, is close to the nonrandom limiting d.f. guaranteed by
Theorem 1 where $m=p+q$, $y=(p+q)/n$, and $H=F^{[B \:\:\:\sigma
I_p]^*[B \:\:\:\sigma I_p]}$.  For the purpose of removing the $q$
singularities, note that the limiting d.f. in Theorem 1 does not
depend on the d.f. of $Y_{11}$. Therefore, without loss of
generality, we may assume that the entries of $V$ are standardized
complex Gaussian. Now, if we let $O^*\Lambda O$ denote the spectral
decomposition of $[B \:\:\: \sigma I_p]^*[B \:\:\:\sigma I_p]$,
where 
the eigenvalues are arranged in nonincreasing order along the
diagonal of $\Lambda$, then the eigenvalues of the matrix in
(\ref{eig5}) are the same as those of
\begin{equation}
\label{those}
\frac1n\Lambda^{1/2}OVV^*O^*\Lambda^{1/2}.
\end{equation}
Since the entries of $V$ are i.i.d. standardized complex Gaussian,
so are the entries of $OV$. Notice the entries of the matrix in
(\ref{those}) outside the upper left $p\times p$ submatrix are
zero. Therefore, the  spectrum of $\widehat R$ is the same as that
of the $p \times p$ upper block. Note also that the $p$ largest
eigenvalues of $[B \:\:\:\sigma I_p]^*[B \:\:\:\sigma I_p]$ are the
same as those of
\begin{equation}
[B\:\:\:\sigma I_p][B\:\:\:\sigma
I_p]^*=BB^*+\sigma^2I_p. 
\end{equation}
Let $P^*\Lambda^{\prime} P$ denote the
spectral decomposition of $BB^*+\sigma^2I_p$, and let $Z$ denote
the first $p$ rows of $OV$. Then, the spectrum of $\widehat R$ is
the same as that of
\begin{equation}
\dfrac{1}{n}ZZ^*\Lambda^{\prime}=\dfrac{1}{n}
ZZ^*P\big(BB^*+\sigma^2I_p\big)P^*, 
\end{equation}
which is the same as that of 
\begin{equation}
\frac{1}{n}Y_pY_p^*(BB^*+\sigma^2I_p)={\widehat R^{\prime}}, 
\end{equation}
where $Y_p=P^*Z$ is $p\times n$ and contains i.i.d.
standardized complex Gaussian entries.

\end{document}